\documentstyle{amsppt}
\magnification 1200
\input pictex
\def\today
{\ifcase\month\or
     January\or February\or March\or April\or May\or June\or
     July\or August\or September\or October\or November\or December\fi
     \space\number\day, \number\year}
\magnification 1200
\UseAMSsymbols
\hsize 5.5 true in
\vsize 8.5 true in
\parskip=\medskipamount
\NoBlackBoxes

\def\mathbb{\Bbb}

\def\mathcal{\Cal}

\def\vp{\varphi}
\def\arrowk{^\to{\kern -6pt\topsmash k}}
\def\arrowK{^{^\to}{\kern -9pt\topsmash K}}
\def\arrowr{^\to{\kern-6pt\topsmash r}}
\def\bark{\bar{\kern-0pt\topsmash k}}
\def\arrowvp{^\to{\kern -8pt\topsmash\vp}}
\def\arrowf{^{^\to}{\kern -8pt f}}
\def\arrowg{^{^\to}{\kern -8pt g}}
\def\arrowu{^{^\to}a{\kern-8pt u}}
\def\arrowt{^{^\to}{\kern -6pt t}}
\def\arrowe{^{^\to}{\kern -6pt e}}
\def\tk{\tilde{\kern 1 pt\topsmash k}}
\def\barm{\bar{\kern-.2pt\bar m}}
\def\barN{\bar{\kern-1pt\bar N}}
\def\barA{\, \bar{\kern-3pt \bar A}}

\def\mathbb{\Bbb}

\def\snint{\raise2pt\hbox{$_{^\not}$}\kern-3.5 pt\int}

\TagsOnRight
\NoRunningHeads

\document
\topmatter
\title
Disjointness of Mobius from horocycle flows
\endtitle

\abstract
  We formulate and prove a finite version of
Vinogradov's bilinear sum inequality. We use
it together with Ratner's joinings theorems to
prove that the Mobius function is disjoint from 
discrete horocycle flows on $\Gamma \backslash SL_2(\Bbb R)$.
\endabstract

\author
J.~Bourgain, P.~Sarnak and T.~Ziegler
\endauthor
\address
School of Math., Institute for Advanced Study, Princeton, NJ 08540, USA \
\endaddress
\email
bourgain\@ias.edu and sarnak\@ias.edu
\endemail
\address
Math. Dept., Technion, Haifa 3200, Israel
\endaddress
\email
tamarzr\@tx.technion.ac.il 
\endemail
\endtopmatter

{\bf 1. Introduction}

In this note we establish a new case of the disjointness conjecture [Sa1] concerning the Mobius function $\mu(n)$.
The conjecture asserts that for any deterministic topological dynamical system $(X, T)$ (that is a compact metric space $X$ with a continuous map
$T$ of zero entropy) as $N\to \infty$,
$$
\sum_{n\leq N} \mu(n) f(T^n x)=o(N)\tag 1.1
$$
where $x\in X$ and $f\in C(X)$.

If this holds we say that $\mu$ is disjoint from $(X, T)$.
The conjecture is known for some simple deterministic systems.
For $(X, T)$ a Kronecker flow (that is a translation in a compact abelian group) it is proven in [V] and [D] while for $(X, T)$ a translation on
a compact nilmanifold it is proved in [G-T].
It is also known for some substitution dynamics associated with the Morse sequence [M-R].\footnote"{$^{(*)}$}"{Also related to this last case is the 
orthogonality of $\mu$ to $AC_0$
functions, see [K], [G], [B].}
In all of these the dynamics is very structured, for example it is not mixing.
Our aim is to establish the conjecture for horocycle flows for which the dynamics is much more random being mixing of all orders $[M]$.

In more detail let $G=SL_2(\Bbb R)$ and $\Gamma\leq G$ a lattice, that is a discrete subgroup of $G$ for which $\Gamma \backslash G$ has finite volume.
Let $u= \left[\matrix 1&1\\ 0&1\endmatrix\right]$ be the standard unipotent element in $G$ and consider the discrete horocycle flow $(X, T)$,
where $X=\Gamma\backslash G$ and $T$ is given by
$$
T(\Gamma x) =\Gamma x u.\tag 1.2
$$

\proclaim
{Theorem 1}
Let $(X, T)$ be a horocycle flow, then $\mu$ is disjoint from $(X, T)$, that is given $x\in X$ and $f\in C(X)$ (if $X$ is not compact then
$f$ is continuous on the one point compactification of $X$), as $N\to \infty$
$$
\sum_{n\leq N}\mu(n) f(T^n x)=o(N).
$$
\endproclaim

\noindent
{\bf Note 1.} We offer no rate in this $o(N)$ statement.
For this reason we cannot say anything about the corresponding sum over primes\footnote"{$^{(*)}$}"{see [S-U] for some results on sums on primes in this case.}
 (which the treatments in the cases of Kronecker
and nilflows certainly do).
The source of the lack of a rate is that we appeal to Ratner's Theorem [R1] concerning joinings of horocycle flows and her proof yields
 no rates.

As pointed out in [Sa1], Vinogradov's bilinear method for studying sums, over primes or correlations with $\mu(n)$, has a natural dynamical interpretation in the context of
 the sequences $f(T^n x)$ belonging to flows.
That is the so called type I sums [Va] are individual Birkhoff sums for $(X, T^{d_1})$ and the type II sums are such Birkhoff sums for joinings of
$(X, T^{d_1})$ with $(X, T^{d_2})$.
The standard treatments ([Va], [I-K]) assume that one has at least a $(\log N)^{-A}$ rate for those dynamical sums in setting up the sieving process.
Our starting point is to formulate a finite version of the bilinear sums method.
It applies to any multiplicative function bounded by 1.

\proclaim
{Theorem 2} 
Let $F:\Bbb N\to\Bbb C$ with $|F|\leq 1$ and let $\nu$ be a multiplicative function with $|\nu|\leq 1$.
Let $\tau>0$ be a small parameter and assume that for all primes $p_1, p_2\leq e^{1/\tau}$, $p_1\not=p_2$, we have that for $M$ large enough
$$
\Big|\sum_{m\leq M} F(p_1 m)\overline{F(p_2m)}\Big| \leq\tau M.\tag 1.3
$$
Then for $N$ large enough
$$
\Big|\sum_{n\leq N}\nu(n)F(n)\Big|\leq 2\sqrt{\tau\log 1/\tau} N.
\tag 1.4
$$
\endproclaim

\noindent
{\bf Note 2.}
There are obvious variations and extensions which allow a small set of $p_1, p_2$ for which (1.3) fails, but for which the conclusion (1.4) may still be
drawn.
We will note them as they arise.

Theorem 2 can be applied to flows $(X, T)$ with $F(n) =f(T^n x)$ as long as we can analyze the bilinear sums $f(T^{p_1n} x) f(T^{p_2n} x)$.
In Section 3 and 4 we use Ratner's theory of joinings of horocycle flows to compute the correlation limits
$$
\lim_{N\to\infty}\frac 1N \sum_{n\leq N} f(T^{p_1 n} x) f(T^{p_2 n} x)
$$
when $(X, T)$ is a horocycle flow.
This correlation is determined by a subgroup of $\Bbb R^*_{>0}$ denoted by $C(\Gamma,  x)$ which is defined in terms of the point $x\in X$ and the
commensurator,  $COM(\Gamma)$ of $\Gamma$ in $G$ (see Section 3).
After removing the mean of $f$ (with respect to $dg$ on $X$) and determining the correlation limits in (1.4) we find that for $T$-generic $x\in X$ (1.3) holds
for $\tau$ as small as we please except for a limited number of $p_1, p_2$'s.
This leads to Theorem 1 if $x$ is generic, while if it is not so then Theorem 1 follows from the Kronecker case.

The method used to handle $G= SL_2(\Bbb R)$ has the potential to apply to the general $Ad$-unipotent flow in $\Gamma\backslash G$, with $G$ semisimple and $\Gamma$ such a
lattice.
For these the correlations are still very structured by Ratner's general rigidity Theorem [R2].
However the possibilities for the correlations are more complicated and we have not examined them in detail.
There are other deterministic flows for which we can apply Theorem 2 such as various substitution flows ([F]) and rank one systems ([Fe]).
We comment briefly on this at the end of the paper, leaving details and work in progress for a future note.

\bigskip

\noindent
{\bf  Section 2. A finite version of Vinogradov's inequality} 

We prove Theorem 2.
The basic idea is to decompose the set of integers in the interval $[1, N]$ into a fixed number of pieces depending on the small parameter $\tau$.
These are chosen to cover most of the interval and so that the members of the pieces have unique prime factors in suitable dyadic intervals.
In this way one can use the multiplicativity of $\nu$ and after an application of Cauchy's inequality, one can invoke (1.3) in order to
estimate the key sum in the Theorem.

Let $\alpha>0$ (small and to be chosen later to depend on the parameter $\tau$) and set
$$
j_0= \frac 1\alpha \Big(\log \frac 1\alpha\Big)^3, j_1= j_0^2, D_0=(1+\alpha)^{j_0} \text { and } D_1=(1+\alpha)^{j_1}.\tag 2.1
$$
In order to decompose $[1, N]$ suitably, consider first the set $S$ given as
$$
S=\{n\in [1, N]: n \text { has a prime factor in $(D_0, D_1)$}\}.\tag 2.2
$$
In what follows $N\to \infty$ with our fixed small $\alpha$ and $A\lesssim B$ means that asymptotically as $N\to\infty$, $A\leq B$.
From the Chinese remainder Theorem it follows that (here and in what follows $[1, N]$ consists of the integers in this interval)
$$
|[1, N]\backslash S|\lesssim \prod_{\Sb D_0<\ell<D_1\\ \ell \text { prime}\endSb} \Big(1-\frac 1\ell\Big)N.\tag 2.3
$$
We can estimate the product over primes in (2.3) using the prime number theorem and that $\alpha$ is small and hence $D_0$ large;
$$
\prod_{D_0<\ell<D_1} \Big(1-\frac 1\ell\Big) \sim \frac{\log D_0}{\log D_1}=\frac 1{j_0}.\tag 2.4
$$
It follows that
$$
|[1, N)\backslash S| \lesssim \alpha N, \tag 2.5
$$
that is up to a fraction of $\alpha, S$ covers $[1, N)$.

Let $P_j$ be the set of primes in $[(1+\alpha)^j, (1+\alpha)^{j+1}]$ for $j_0\leq j\leq j_1$ and define $S_j$ by
$$
S_j=\{n\in [1, N); n\text { has a single divisor in $P_j$ and no divisor in $\bigcup_{i<j}P_i$}\}.\tag 2.6
$$
The sets $S_j$ are disjoint and appealing again to the prime number theorem with remainder and $\alpha$ small we have
$$
|P_j|=\frac {(1+\alpha)^{j+1}}{(j+1)\log (1+\alpha)} -\frac{(1+\alpha)^j}{j\log(1+\alpha)}+ O\big((1+\alpha)^j e^{-\sqrt{\alpha j}}\big),
\tag 2.7
$$
with an implied constant that is absolute.

Hence for $\alpha$ small and $j_0\leq j\leq j_1$,
$$
|P_j|\leq (1+\alpha)^j \Big[\frac 1j+\frac 1{\alpha j^2}+O(e^{-\sqrt{\alpha j}})\Big].\tag 2.8
$$
Now from the definition of $S$ we have that
$$
S\backslash\bigcup_{j_0\leq j\leq j_1} S_j \subset\bigcup_{j_0\leq j\leq j_1} \{n\in [1, N); \ 
\matrix {}\\ {\text{with  $n$ having at least }}\\
{\text{ \ two prime factors in $S_j$}\}}.\endmatrix
\tag 2.9
$$
Hence
$$
\align
|S\backslash\bigcup_{j_0\leq j\leq j_1} S_j| &\lesssim \sum_{j_0\leq j\le j_1} \ \sum_{\ell_1, \ell_2 \in P_j} \frac N{\ell_1\ell_2}
\\
&\leq N\sum_{j_0 \leq j\leq j_1} \Big(\frac {|P_j|}{(1+\alpha)^j}\Big)^2.\tag 2.10
\endalign
$$
From (2.8) this gives for $\alpha$ small enough
$$
\align
|S\backslash \bigcup_{j_0\leq j\leq j_1} S_j|&\lesssim N\sum_{j_0\leq j\leq j_1} \Big(\frac 1j+\frac 1{\alpha j^2}+O\big(e^{-\sqrt {\alpha j}}\big)
\Big)^2\\
&\leq N\Big(\frac 1{j_0} +\frac 1{j_0^3\alpha^2}+ O\Big(\frac 1\alpha (1+\sqrt{\alpha j_0}) e^{-\sqrt{\alpha j_0}}\Big)\Big)\\
&\leq \alpha N.\tag 2.10$'$
\endalign
$$

So at this point we have covered $[1, N]$ up to a fraction of $\alpha$ by the disjoint sets $S_j, j_0\leq j\leq j_1$.
Finally we decompose $S_j$ in a well factored set and its complement.
For $j_0\leq j\leq j_1$, let
$$
Q_j=\Big\{m\in \Big[1, \frac N{(1+\alpha)^{j+1}}\Big); \text { $m$ has no prime factors in $\bigcup_{i\leq j}P_j$}\Big\}.\tag 2.10$''$
$$
Clearly the product sets $P_jQ_j$ satisfy
$$
P_jQ_j\subset S_j \text { for } j_0\leq j\leq j_1.\tag 2.11
$$
Moreover for each such $j$
$$
S_j\backslash (P_jQ_j) \subset P_j .
\Big[\frac N{(1+\alpha)^{j+1}}, \frac N{(1+\alpha)^{j}}\Big].\tag 2.12
$$
Hence
$$
\sum_{j_0\leq j\leq j_1} |S_j\backslash (P_jQ_j)|\leq \sum_{j_0\leq j\leq j_1} |P_j|.\frac {\alpha N}{(1+\alpha)^j}.
$$
Applying (2.7) yields that the right hand side above is
$$
\leq N\Big\{\alpha \log \frac {j_1}{j_0}+\frac 1{j_0}+O(1+\sqrt{\alpha j_0} e^{-\sqrt{\alpha j_0}})\Big\}.
$$
Hence for $\alpha$ small enough
$$
\sum_{j_0\leq j\leq j_1} |S_j\backslash (P_j Q_j)|\leq 2\alpha N.\tag 2.13
$$
This leads to the basic decomposition of $[1, N)$ into disjoint sets $P_jQ_j, j_0\leq j\leq j_1$ with only a small number of
points omitted.
Namely from (2.10$''$) and (2.13)
$$
|[1, N)\backslash\bigcup_{j_0\leq j\leq j_1} P_jQ_j|\lesssim 3\alpha N.\tag 2.14
$$
Note that the map $P_j\times Q_j \to P_j Q_j$ is one-to-one and since $|F|\leq 1$ and $|\nu|\leq 1$ we have that
$$
\Big| \sum_{n\leq N} \nu(n) F(n)\Big| \lesssim \sum_{j_0\leq j\leq j_1} \Big| \sum_{\Sb x\in P_j\\ y\in Q_j\endSb} \nu(xy)F(xy)\Big|+3\alpha N.
\tag 2.15
$$
For $x\in P_j, y\in Q_j, (x, y)=1$ so that the $\nu(xy)$ in (2.15) can be factored as $\nu(x)\nu(y)$ and hence
$$
\Big| \sum_{\nu\leq N} \nu(n)F(n)\Big|\lesssim \sum_{j_0\leq j\leq j_1} \ \sum_{y\in Q_j}\Big|\sum_{x\in P_j} \nu(x) F(xy)\Big|+3\alpha N.\tag 2.16
$$
The inner sum may be estimated using Cauchy;
$$
\align
&\sum_{y\in Q_j} \Big|\sum_{x\in P_j} \nu (x) F(xy)\Big|\\
&\leq \Big(\sum_{y\in Q_j}  1\Big)^{1/2} \Big(\sum_{y\in Q_j} \Big|\sum_{x\in P_j} \nu(x)F(xy)\Big|^2\Big)^{1/2}\\
&\leq |Q_j|^{1/2} \Big(\sum_{y\leq N/(1+\alpha)^j}\Big|\sum_{x\in P_j} \nu(x) F(xy)\Big|^2\Big)^{1/2}\\
&=|Q_j|^{1/2}\Big(\sum_{y\leq N/(1+\alpha)^j} \  \sum_{x_1, x_2\in P_j} \nu(x_1) \overline{\nu(x_2)} F(x_1 y) \overline{F(x_2 y)}\Big)^{1/2}\\
&\leq |Q_j|^{1/2}\Big(\sum_{x_1, x_2 \in P_j}\Big| \sum_{y\leq N/(1+\alpha)^j} F(x_1 y) \overline{F(x_2y)}\Big|\Big)^{1/2},\tag 2.17
\endalign
$$
where we have used $|\nu|\leq 1$.

Note that here
$$
x_1, x_2 <(1+\alpha)^{j_1}< e^{1/\alpha^2}.\tag 2.18
$$
The diagonal contribution in (2.17), that is $x_1=x_2$ for each $j$ yields at most
$$
|Q_j|^{1/2}|P_j|^{1/2} \frac {\sqrt N}{(1+\alpha)^{j/2}},\tag 2.19
$$
by using that $|F|\leq 1$ and the definition of $Q_j$.
Hence summing over $j$ and Cauchy it is 
$$
\leq \sqrt N\Big(\sum_{j_0\leq j\leq j_1} |P_j| \, |Q_j|\Big)^{1/2} \Big(\sum_{j_0\leq j\leq j_1} \frac 1{(1+\alpha)^j}\Big)^{1/2}.
$$
Now $|P_jQ_j|\leq |S_j|$ and $\sum|S_j|\leq N$, hence the full diagonal contribution is at most
$$
N\Big(\sum_{j_0\leq j\leq j} \frac 1{(1+\alpha)^j}\Big)^{1/2} \leq \alpha N.\tag 2.20
$$
For $x_1\not= x_2$ the hypothesis in the Theorem may be applied in view of (2.18), that is
$$
\Big|\sum_{y\leq N/(1+\alpha)^j} F(x_1 y) \overline{F(x_2y)}\Big|\leq \frac {\tau N}{(1+\alpha)^j}.
$$

Hence the off-diagonal contribution is at most
$$
\align
&\sqrt{\tau N}\sum_{j_0\leq j\leq j_1} |P_j| \, |Q_j|^{1/2} (1+\alpha)^{-j/2}\\
&\leq\sqrt{\tau N} \Big(\sum_{j_0\leq j\leq j_1} |P_j|\, |Q_j|^{1/2} \Big)^{1/2} \Big(\sum_{j_0\leq j\leq j} |P_j|(1+\alpha)^{-j}\Big)^{1/2}\\
&\leq \sqrt{\tau N} N^{1/2} \Big(\log \frac {j_1}{j_0}+\frac 1{j_0\alpha} +\frac 1\alpha (1+\sqrt{\alpha j_0}) e^{-\sqrt{j_0}}\Big)^{1/2}\\
&\le N\sqrt\tau\sqrt{\log 1/\alpha}\tag 2.21
\endalign
$$
(for $\alpha$ small).

Putting all of these together we have
$$
\Big|\sum_{n\leq N} \nu(n) F(n)\Big|\lesssim N\big(4\alpha+\sqrt{\tau\log 1/\alpha}\big).
$$
Taking $\alpha=\sqrt\tau$ yields the Theorem.

\bigskip
\noindent
{\bf Section 3. Commensurators and correlators}

As in the introduction $G=SL_2(\Bbb R)$ and $\Gamma$ is a lattice in $G$.
The commensurator subgroup, $COM(\Gamma)$ of $\Gamma$ in $G$, is defined by
$$
COM(\Gamma) =\{ g\in G: g^{-1} \Gamma g\cap\Gamma \text { is finite index in both $\Gamma$ and $g^{-1}\Gamma g$}\}.\tag 3.1
$$
It plays a critical role in determining the ergodic joinings of $(\Gamma\backslash G, T^a)$ with $(\Gamma\backslash G, T^b)$, where
$a, b>0$ and $T^a=\left[\matrix 1&a\\ 0&1\endmatrix\right]$.
Let $z$ be a point on the projective line $\Bbb P^1(\Bbb R)$ and let $P_z$ be the stabilizer of $z$ in $G$, with $G$ acting projectively.
If $z=\infty$ then
$$
P_\infty=\Big\{\pmatrix \alpha&\beta\\ 0&\delta\endpmatrix:\alpha \delta =1\Big\}.\tag 3.2
$$
If $\xi\in G$ and $\xi(\alpha)=z$ then
$$
P_z =\xi^{-1} P_\infty \xi.\tag 3.3
$$
Define the character $\chi$ of $P_\infty$, and hence of $P_z$ for any $z$, by
$$
\chi\Big(\pmatrix \alpha&\beta\\0&\delta\endpmatrix\Big) =\alpha\delta^{-1} =\alpha^2.\tag 3.4
$$
$\chi$ is valued in the multiplicative group $\Bbb R^*_{>0}$.
If $\Delta_z$ is a subgroup of $P_z$ we define the correlation group $C(\Delta_z)$ to be the
image of $\Delta_z$ under $\chi_z$, that is $C(\Delta_z)$ is the subgroup of $\Bbb R^*$ given as $\chi_z(\Delta_z)$.
We denote by $C(\Gamma, z)$ the group $C\big((COM \, \Gamma)\cap P_z)$ and our aim is to determine this group for $\Gamma$
and $z$ as above.
Its relevance to the unipotent element $u$ is that for $\beta\in P_\infty$
$$
\beta u\beta^{-1}=\left[\matrix 1& \chi(\beta)\\ 0&1\endmatrix\right] = u^{\chi(\beta)}.\tag 3.5
$$

The explicit computation of these groups $C(\Gamma, z)$ depends on the nature of $\Gamma$, so we divide it into cases.

\noindent
{\bf Case 1.} In which $\Gamma$ is nonarithmetic and $\Gamma\backslash G$ is compact.
In this case it is known [Ma] that $COM(\Gamma)/\Gamma$ is finite and hence $COM(\Gamma)$ is itself a lattice in $G$.
Hence for $z\in\Bbb P^1 (\Bbb R), COM(\Gamma)\cap P_z$ is cyclic (either trivial or infinite) and hence what is important
for us is that $C(\Gamma, z)$ is finitely generated.  
In particular it follows that the set of $p/q$ with $p\not= q$ and both prime which lie in $C(\Gamma, z)$ is finite.
We record this as

\proclaim
{Lemma 1} If \, $\Gamma\backslash G$ is compact and $\Gamma$ is nonarithmetic then for any $z\in\Bbb P^1(\Bbb R)$
$$
\Big\{ \frac pq: p, q \text { prime } p\not= q\Big\} \bigcap C (\Gamma, z)
$$
is finite (in fact consists of at most one element).
\endproclaim

\noindent
{\bf Case 2.}
In which $\Gamma\backslash G$ is compact and $\Gamma$ is arithmetic.
In this case it is known [We] that $\Gamma$ is commensurable with a unit group in a quaternion algebra $A$ defined over a totally real number
field.
For simplicity we assume that $A$ is defined over $\Bbb Q$ (the general case may be analyzed similarly).
Thus
$ A=\Big(\frac {a, b}Q\Big)$ is a 4-dimensional division
algebra (since $\Gamma\backslash G$ is compact) generated linearly over $\Bbb Q$ by $1, \omega, \Omega, \omega\Omega$.
Here $\omega^2 =a, \Omega^2=b$ with $a, b\in \Bbb Q$ and say $a>0$ and $a$ and $b$ square-free.
$\omega$ and $\Omega$ obey the usual quaternionic multiplication rules.
For
$$
\alpha =x_0+x_1 \omega+x_2\Omega+x_3 \omega\Omega\tag 3.6
$$
with $x_j\in\Bbb Q$
$$
\overline\alpha = x_0-x_1 \omega-x_2\Omega -x_3 \omega\Omega,\tag 3.7
$$
$$
N(\alpha)=\alpha \overline\alpha =x_0^2- ax_1^2 -bx_2^2 +abx_3^2\tag 3.8
$$
and
$$
\text {trace $(\alpha) =\alpha +\overline \alpha =2x_0$}.\tag 3.9
$$
$A/\Bbb Q$ being a division algebra  is equivalent to the statement
$$
N(\alpha)=0 \text { iff  $\alpha=0$, for $\alpha \in A(\Bbb Q)$}.\tag 3.10
$$
Let $A_1(\Bbb Z)$ be the integral unit group;
$$
A_1(\Bbb Z) =\{\alpha\in A(\Bbb Z):N(\alpha)=1\}.\tag 3.11
$$
We embed $A(\Bbb Q)$ into $M_2(\Bbb R)$ by
$$
\alpha \to \phi(\alpha) =\left[ \matrix {\overline\xi}& \eta\\ b {\overline \eta} &\xi\endmatrix\right] \tag 3.12
$$
where $\xi =x_0-x_1 w, \eta= x_2+x_3 \omega$ and $\omega=\sqrt a\in\Bbb R$. 

Note that
$$
\det \phi(\alpha) =N(\alpha)\tag 3.13
$$
$$
\text {trace\,} \phi(\alpha)=\xi+\overline\xi = \text { trace $(\alpha)$}.\tag 3.14
$$
Now $\Lambda =\phi\big(A_1(\Bbb Z)\big)$ is a cocompact lattice in $G$ and we are assuming that our  $\Gamma$ is commensurable with $\Lambda$.
Hence the commensurator of $\Gamma$ (or of $\Lambda$ they are the same) consists of the $\Bbb Q$ points [P-R];
$$
COM (\Gamma) =\left\{ \frac {\phi(\alpha)}{(\det \alpha)^{1/2}}; \alpha\in A^+(\Bbb Q)\right\}\tag 3.15
$$
where $A^+(\Bbb Q)$ consists of all $\alpha\in A(\Bbb Q)$ with $N(\alpha)>0$.

Hence up to scalar multiples of $\left[\matrix 1&\\ &1\endmatrix\right]$, that is up to the center of $GL_2(\Bbb R)$
$$
\delta\in COM (\Gamma) \text { iff } \delta =\left[\matrix \overline\xi&\eta\\ b\overline\eta&\xi\endmatrix\right]
\text { with  $\xi +\eta\Omega\in A^+(\Bbb Q)$.}\tag 3.16
$$

Our interest is in $C(\Gamma, z)$ and from the description (3.16) one can check that for certain algebraic $z$'s this group can be an infinitely
generated subgroup of $K^*$, where $K$ is the corresponding algebraic extension of $\Bbb Q$.

What is important to us are the rationals in this group and this is given by

\proclaim
{Lemma 2}
For $\Gamma$ as in case (ii) and any $z\in\Bbb P^1(\Bbb R)$
$$
C(\Gamma, z) \bigcap \Bbb Q^* =\{1\}.
$$
\endproclaim

\noindent
{\bf Proof.} Let $\hat\delta \in P_z\bigcap COM (\Gamma)$, then $\hat\delta =\phi(\delta)$ and hence $N(\delta)^{1/2} \hat\delta$ in $GL_2(\Bbb R)$
satisfies
$$
\left.
\aligned
\text {trace } (N(\delta)^{1/2} \hat\delta) &= s\in\Bbb Q\\
\det (N(\delta)^{1/2} \hat\delta)&= t\in\Bbb Q^*_{>0}\endaligned
\right\}
\tag 3.17
$$
Also $N(\delta)^{1/2} \hat\delta$ is conjugate in $G$ to $\beta$
with
$$
\beta=\left[\matrix \lambda&*\\ 0&\mu\endmatrix\right]
$$
where
$$
\lambda\mu=t \text { and } \lambda+\mu=2s.\tag 3.18
$$
Now $\chi(\hat\delta)=\lambda/\mu$ and if this number is in $\Bbb Q$ then from (3.17) and (3.18) we see that both $\lambda$ and $\mu$ are
in $\Bbb Q$.
Now $\delta\in A^+(\Bbb Q)$ so $\delta=x_0+x, \omega+x_2\Omega +x_3\omega\Omega$ with $x_j\in \Bbb Q$ and from (3.17) we have that
$$
\Big(\frac {\lambda+\mu} 2\Big)^2 -ax_1^2 -bx_2^2 +abx_3^2=\lambda\mu
$$
that is
$$
\Big(\frac {\lambda-\mu}2\Big)^2 -a_1 x_1^2 -bx_2^2 +abx_3^2 =0.\tag 3.19
$$
Now $\lambda-\mu, x_1, x_2, x_3 \in\Bbb Q$ and  since $A$ is a division algebra it follows from (3.19) that
$$
\lambda -\mu =x_1= x_2=x_3=0.
$$
That is $\lambda =\mu$ and hence $\lambda/\mu=1$ as claimed.

\noindent
{\bf Case 3.} In which $\Gamma\backslash G$ is a finite volume, but noncompact and nonarithmetic.
As in the cocompact nonarithmetic case we have that $COM(\Gamma)$ is a lattice in $G$.
Hence as before $P_z\cap COM(\Gamma)$ is cyclic and in particular $C(\Gamma, z)$ is a finitely generated subgroup of $\Bbb R^*>0$
and Lemma 1 is true for such a $\Gamma$.
This leaves us with the last case.

\noindent
{\bf Case 4.} In which $\Gamma$ is arithmetic and noncompact.
This time $\Gamma$ is commensurable with a quaternion algebra that is split over $\Bbb Q$ and hence $\Gamma$ is commensurable with
$SL_2(\Bbb Z)$.
Its commensurator subgroup is given by
$$
COM (\Gamma)=\{A/(\det A)^{1/2}: A\in GL_2^+(\Bbb Q)\}.
$$
Now if $z\in\Bbb P^1(\Bbb Q)$, then since
$$
COM (\Gamma) \cap P_\infty =\Big\{ \frac 1{\sqrt{\alpha\delta}} \pmatrix \alpha&\beta\\ 0\, &\delta\endpmatrix: \alpha, \beta, \delta\in\Bbb Q,
\alpha\delta>0\Big\}
$$
we have that
$$
C(\Gamma, z) =C \big(COM(\Gamma)\cap P_z\big) =\Bbb Q^* .\tag 3.20
$$
So in this case the correlator subgroup contains every rational $p/q$.

If $z\not\in \Bbb P^1 (\Bbb Q)$ and $z$ does not lie in a quadratic number field then $(az+b)/ (cz+d)=z$ has no solutions for $\left[\matrix a&b\\
c&d\endmatrix\right] \in GL(2, \Bbb Q)$ other than $\left[\matrix a&b\\c&d\endmatrix\right] =1\!\!1$ in $PGL_2$.
Hence for such a $z$
$$
C(\Gamma, z) =C\big(COM (\Gamma)\cap P_z\big)=\{1\}.\tag 3.21
$$
This leaves us with $z$ quadratic in which case we have that up to scalar multiples of $I$;
$$
COM(\Gamma)\cap P_z=\{\gamma\in GL_2^+(\Bbb Q):\gamma z=z\}.\tag 3.22
$$
If $z$ satisfies $az^2+bz+c=0$ with $a, b, c$ integers $(a, b, c)=1$ and $d=b^2-4ac>0$
and not $a$ square, then one checks that
$$
COM(\Gamma)\cap P_z =\left\{\pmatrix \frac {t+ba} 2& au\\ -cu& \frac {t-bu} 2\endpmatrix : t^2-du^2\in\Bbb Q^+  \, t, u\in \Bbb Q\right\}.
$$

Hence
$$
\chi\Bigg(\pmatrix \frac {t-bu} 2 & au\\ -cu&\frac {t+bu}2\endpmatrix\Bigg) =\frac {t+u\sqrt d}{t-u\sqrt d}\tag 3.23
$$
and so
$$
C(\Gamma, z)=\left\{\frac \eta {\eta'}; \eta\in \Bbb Q(\sqrt d)^*, N(\eta) >0\right\}\tag 3.24
$$
where $\eta'$ is the conjugate of $\eta$ in $\Bbb Q(\sqrt d)$.

While this groups is infinitely generated its intersection with $\Bbb Q^*$ is 1
(if $(\frac \eta{\eta'})' =\frac \eta{\eta'}$.
Then $\mu^2 =(\eta')^2 $ or $\eta=\pm\eta'$ and since $N(\eta)>0, \eta=\eta'$).

We summarize this with

\proclaim
{Lemma 3}  If $\Gamma$ is commensurable with $SL_2(\Bbb Z)$ then
$$
C(\Gamma, z) \cap\Bbb Q^*= \left\{\aligned & \{1\} \text { if } z\not\in \Bbb P^1(
\Bbb Q)\\
&\Bbb Q^* \text { if } \ z\in  \Bbb P^1(\Bbb Q).\qquad\qquad
\endaligned\right.
$$
\endproclaim

\bigskip
\noindent
{\bf Section 4. Ratner rigidity and Mobius disjointness}

The correlator group $C(\Gamma, z)$ enters in the analysis of joinings of horocycle flows when applying Ratner's Theorem ([R1], [R2]).
According to these we have that for $\lambda_1, \lambda_2>0$ and $\xi\in\Gamma\backslash G$
$$
\lim_{N\to\infty} \frac 1N\sum^N_{n=1} f\big(\xi(u^{\lambda_1})^n\big) f \big(\xi(u^{\lambda_2})^n\big)\tag 4.1
$$
exists.
Here $f\in C(\Gamma\backslash G)$ is continuous on the one point compactification of $\Gamma\backslash G$ (if it is not compact) and $u^\lambda =\left
[\matrix 1&\lambda\\ 0&1\endmatrix \right]$.

The limit in (4.1) is given by
$$
\int_{\Gamma\backslash G\times \Gamma\backslash G} F(\tilde\xi h) d\nu (h)\tag 4.2
$$
where $\tilde\xi =(\xi, \xi)\in X\times X, F(x_1, x_2)= f(x_1) f(x_2)$ and $\nu$ is an algebraic Haar measure supported on an algebraic
subgroup $H$ of $G\times G$ and for which $(\Gamma\times\Gamma)\tilde \xi H$ is closed in $X\times X$.
The support of $\nu$ is the closure of the orbit $\big(\xi(u^{\lambda_1})^n, \xi(u^{\lambda_2})^n\big), n=1, 2, \ldots$.

Consider first the case that $X$ is compact.
Then any point $\xi$ is $u^\lambda$ generic (the flow $(X, T^\lambda)$ is uniquely ergodic and every point is $dg$ equidistributed in $X$ [Fu]).
It follows that the measure $\nu$ projects onto $dg$ on each factor $X_{\lambda_j}$ of $X_{\lambda_1}\times X_{\lambda_2}$.
That is $\nu$ is a joining of $X_{\lambda_1}$ with $X_{\lambda_2}$.
Applying the Ratner rigidity Theorems either $d\nu =dg_1\times dg_2$ or it reduces to a measure on subgroups $H=\psi(G)$ where $\psi$ is a morphism
$\psi:G \to G\times G$  of the form
$$
\psi(g)=\big(\psi_1 (g), \psi_2(g)\big)
$$
with
$$
\psi_1 (u) =u^{\lambda_1}, \psi_2(u) =u^{\lambda_2}\tag 4.3
$$
and $\big(\Gamma \xi\psi_1 (g), \Gamma\xi\psi_2(g)\big)$ is closed in $X\times X$.
That is there are $\alpha_1, \alpha_2 \in G$ such that
$$
\alpha_1 u  \alpha_1^{-1} =u^{\lambda_1}, \alpha_2 u\alpha_2^{-1}=u^{\lambda_2}
$$
and
$$
\psi(g)=(\alpha_1 g\alpha_1^{-1}, \alpha_2 g\alpha_2^{-1}).\tag 4.4
$$
In particular $\alpha_1, \alpha_2 \in P_\infty(\Bbb R)$ and 
$$
\chi(\alpha_1)=\lambda_1\text { and $\chi(\alpha_2)=\lambda_2$}.\tag 4.5
$$
Now $(\xi\alpha_1 g\alpha_1^{-1}, \xi \alpha_2 g\alpha_2^{-1}); g\in G$ is closed in $\Gamma\backslash G\times\Gamma\backslash G$ iff
$$
(h\alpha_1^{-1}, \xi \alpha_2\alpha_1^{-1} \xi^{-1} h\alpha_2^{-1}), h\in G\tag 4.6
$$
is closed $\Gamma\backslash G\times\Gamma\backslash G$.

The latter is equivalent to
$$
\delta =\xi \alpha_2\alpha_1^{-1} \xi^{-1} \in COM (\Gamma).\tag 4.7
$$
In this case
$$
\delta\in P_{\xi(\alpha)} \cap COM (\Gamma) \text { and } \chi(\delta)=\lambda_2/\lambda_1.
\tag 4.8
$$
Thus we have

\proclaim
{Lemma 4}
There is a nontrivial joining (in particular $\nu$ is not $dg_1 dg_2$) in (4.2) iff $\lambda_2/\lambda_1\in C\big(\Gamma, \xi (\infty)\big)$.
\endproclaim

So if $\lambda_1=p$ and $\lambda_2=q$ with $p\not= q$ primes, then from Lemmas 1, 2 and 4 we have;

\proclaim
{Corollary 5}
For $\xi \in\Gamma\backslash G$ with $\Gamma\backslash G$ compact, there are most finite number of pairs of distinct primes $p, q$ (depending only on $\xi$)
for which the following fails
$$
\frac 1N \sum^N_{n=1} f(\xi u^{pn}) f(\xi u^{q n})\to \Big(\int_{\Gamma\backslash G} f(g) dg\Big)^2.
$$
\endproclaim

If $X$ is noncompact and $\Gamma$ is nonarithmetic then as long as $\xi$ is generic for $(X, T, dg)$,
then the joinings analysis coupled with the Case 3 discussion of Section 3 lead to Corollary 5 holding for such $\Gamma$ and $\xi$.
The remaining case is that of $\Gamma$ being commensurable with $SL_2(\Bbb Z)$.
In this case by [Da], $\xi$ is generic for $(X, T, dg)$ iff $\xi(\infty)\not\in \Bbb P^1 (\Bbb Q)$.

Thus again we can apply the joinings analysis coupled with Lemma 3 and Lemma 4 to conclude

\proclaim
{Corollary 6} If $X$ is noncompact and $\xi$ is generic for $(X, T, dg)$ then there are at most a finite number of pairs of distinct prime
$p, q$, depending only on $\xi$, for which the following fails
$$
\frac 1N \sum^N_{n=1} f(\xi u^{pn}) f(\xi u^{qn})\to \Big(\int_{\Gamma\backslash G} f( g)dg\Big)^2.
$$
\endproclaim

We can now complete the proof of Theorem 1.
If $X$ is compact then every $\xi$ is generic for $dg$.
Write
$$
f(x) =f_1(x)+c\tag 4.9
$$
where $\int_{\Gamma\backslash G} f_1(x) dx=0 $ and $c$ is a constant.
Then
$$
\align
&\frac 1N\sum^N_{n=1} \mu(n) f(\Gamma\xi u^n)\\
&=\frac 1N \sum^N_{n=1} \mu(n) f_1 (\Gamma\xi u^n)+o(1)\tag 4.10
\endalign
$$
by the prime number theorem.

As for as the sum against $f_1$ is concerned, according to Corollary 5
(note $\int_{\Gamma\backslash G} f_1(x) dx=0$), 
the conditions  of Theorem 2 are met for $F(n)=f_1(\xi u^n)$ except for finitely many pairs $p, q$.
This causes no harm as for as concluding that the first sum in (4.10) is $o(1)$.
One can certainly allow a finite number of exceptions (independent of $N$) in Theorem 2, in fact the proof only involves the
condition for primes $p\geq D_0$ which gets large as $\tau$ gets small.

If $\Gamma\backslash G$ is not compact and $\xi$ is generic for $dg$ then according to Corollary 6 everything goes through as
above and Theorem 1 follows.  If $\xi$ is not generic then by [Da] the closure of the orbit of $\xi$ in $X$ is either finite or
is a circle and in the latter case the action of $u$ is by rotation of this circle through an angle of $\theta$.
Thus in the first case Theorem 1 follows from the theory of Dirichlet $L$-functions while in the second case it was proven
in [D].

\noindent
{\bf Note 4.} The case of richest joinings of $X\times X$ of the form
$ (\Gamma g\alpha_1^{-1}, \Gamma\delta g\alpha_2^{-1})$ with $\Gamma = SL_2(\Bbb Z)$ and $\delta\in COM (\Gamma)$, is not one
that we had to consider directly in our analysis (since it corresponds to $\xi(\infty) \in \Bbb P^1(\Bbb Q)$ so that $\xi$ is not generic).
For this joining if $\det \delta= p q$ (taking $\delta\in GL_2^+)$
.the joining is
$$
\frac 1{[\Gamma:\Delta]} \int_{\Delta\backslash G} f(g\alpha_1^{-1})f(\delta g\alpha_1^{-1}) dg\tag 4.11
$$
where $\Delta = \delta^{-1} \Gamma\delta \cap\Gamma$.

By the theory of correspondences (Hecke operators) if $f$ is a Hecke eigenform (and $\int_{\Gamma\backslash G} f dg =0)$ which we
can assume here, the joining in (4.11) becomes
$$
\frac {\lambda_ f(pq)}{(p+1)(q+1)} \int_{\Gamma\backslash G} f(g\alpha_1^{-1}) f(g\alpha_1^{-1}) dg
\tag 4.12
$$
where $\lambda_f(n)$ is the $n$-th Hecke eigenvalue.

So while in this case the correlation need not be zero it is small if $pq$ gets large.
This follows from the well known bounds for Hecke eigenvalues [Sa2].
One would expect that this would be useful in an analysis of this type but apparently for this ineffective
analysis it is not needed.

\bigskip
\noindent
{\bf Section 5. Some further comments}

The Moebius orthogonality criterion provided by Theorem 2 has applications to other systems of zero entropy.
One can use it to give a ``soft'' proof of the qualitative Theorem 1 for Kronecker and nilflows.
In what follows we will only briefly review some new consequences that are essentially immediate from a number of classical
facts in ergodic theory, \footnote"{$^{(*)}$}" {We are grateful to J-P.~Thouvenot for a detailed account of
the `state of the art' of various aspects of the theory of joinings.}
leaving details and further research in this direction for a future paper.
Some unexplained terminology below may be found in [Ka-T].
First, we mention a result due to Del Junco and Rudolph ([D-R], Cor. 6.5) asserting the disjointness of distinct powers $T^m$ and $T^n$ for weakly
mixing transformations $T$ with the minimal self-joinings property (MSJ). 
This provides another general class of systems for which Theorem 1 holds.
More precisely, the disjointness statement
of Theorem 1 applies to any uniquely ergodic topological model for such transformations, these exist by Jewett [J].
Next, restricting ourselves to rank-one transformations, J.~King's theorem [K1] states that mixing rank-one implies MSJ (well-known examples
include the Ornstein rank-one constructions and Smorodinsky-Adams map, see [Fe]).
The condition of mixing may be weakened to `partial mixing', see [K-T2].
While it seems presently unknown whether any mildly mixing rank-one transformation has MSJ, this property was established in certain other 
cases, such as Chacon's transformation [D-R-S] (which is mildly but not partially mixing).
It was shown that `typical' interval exchange transformations are never mixing
([Ka]), rank-one ([Ve2]), uniquely ergodic ([Ve1],[Mas]), and weakly mixing ([A-F]).
Whether they satisfy Theorem 1 is an interesting question, especially in view of the fact that they are the immediate generalization
of circle rotations.

Finally, the Moebius (and Liouville sequence) orthogonality with sequences arising from substitution dynamics has its importance
from the perspective of symbolic complexity (see [Fe2] for a discussion).
Our approach via Theorem 2 is applicable for sequences produced by an `admissible $q$-automation' (see [Quef] for definitions),
provided its spectral type is of intermediate dimension.
The spectral measure is indeed known to be $(xq)$-invariant, and disjointness of $T^{p_1}$ and $T^{p_2}$ for $p_1\not= p_2$ may in this case
be derived from [Ho].


\Refs
\widestnumber\no{XXXXX}

\ref\no{[A-F]} \by A.~Avila, G.~Forni
\paper Weak mixing for interval exchange transformations and translation flows.
\jour Annals 165 (2007), 637--664.
\endref

\ref\no{[B]}\by J.~Bourgain
\jour ``Moebius-Walsh correlation bounds and an estimate of Mauduit and Rivat'', preprint 2011
\endref

\ref\no {[D]}\by H.~Davenport
\paper Quat
\jour J.~Math. 8, 313-320 (1937)
\endref

\ref\no{[Da]} \by S.~Dani
\paper On uniformly distributed orbits of certain horocycle flows.
\jour
Ergodic Theory and Dynamical Systems (1982), 139--158
\endref

\ref\no {[D-R-S]} \by A.~Del Junco, A.~Rahe, L.~Swanson
\paper Chacon's automorphisms has minimal self-joinings
\jour J.~Anal. Math. 37 (1980), 276--284
\endref

\ref\no {[D-R]} \by A.~Del Junco, D.J.~Rudolph
\paper On ergodic actions whose self-joinings are graphs
\jour ETDS 7 (1987), 531--557
\endref


\ref\no{[F]} \by P.~Fogg
\jour Lecture Notes in Math. 1794, (2002)
\endref

\ref\no {[Fe1]}\by S.~Ferenczi
\paper Systems of finite rank
\jour Colloq, Math. 73 (1997), No 1, 35--65
\endref

\ref\no{[Fe2]}
\paper Rank and symbolic complexity
\jour ETDS 16 (1996), No 4, 663--682
\endref

\ref\no{[Fu]} \by H.~Furstenberg
\jour Lecture Notes in Math 318, 95--115 (1973)
\endref

\ref\no{[G]}\by B.~Green
\jour ``On (not) computing the Mobius function using bounded depth circuits'', preprint 2011
\endref

\ref\no{[G-T]}\by B.~Green and T.~Tao
\paper The Mobius Function is strongly orthogonal to nilsequences
\jour to appear in Annals
\endref

\ref\no{[Ho]}\by M.~Hochman
\paper Geometric rigidity of $xm$ invariant measures
\jour preprint
\endref

\ref\no{[I-K]}\by H.~Iwaniec and E.~Kowalski
\jour Analytic Number Theory, AMS Coll 53 (2004)
\endref

\ref\no{[J]} \by R.I.~Jewett\jour J.~Math.~Mech. 19 (1970) 717-729
\endref

\ref\no{[K]}\by G.~Kalai
\jour ``The ACO Prime Number Conjecture'', http://jilkalai.wordpress.com/2011/ 02/21/ the -ac0-prime-number-conjecture/
\endref

\ref\no{[Ka]} \by A.~Katok
\paper Interval exchange transformations and some special flows are not mixing. 
\jour Israel J. Math 35, no.4, 301-310 (1980)
\endref

\ref\no{[Ka-T]} \by A.~Katok and J.P.~Thouvenot
\paper Spectral properties and combinatorial constructions in ergodic theory
\jour Handbook of Dynamical System, Vol. 1B, 649--743 (2006)
\endref

\ref\no{[Ki]}\by J.~King
\paper Joining rank and the structure of finite rank mixing transformations
\jour J.~Anal. Math. 51 (1988), 182--227
\endref


\ref\no{[Ki-T]}
\by J.~King and J-P.~Thouvenot
\paper An canonical structure theorem for finite joining-rank maps
\jour J.~Analyse Math. 56 (1991), 211--230
\endref

\ref\no{[M-R]} \by C.~Mauduit and J.~Rivat
\jour Annals 171 (2010), 1591--1646
\endref

\ref\no{[Ma]} \by G.~Margulis
\jour Discrete Subgroups of Lie Groups, 1991
\endref

\ref
\no{[M]} \by B.~Marcus
\jour Invent. 46, 201-209 (1978)
\endref

\ref
\no[Mas] \by H.~Masur
\paper Interval exchange transformations and measured foliations.
\jour Annals Vol 115, 210--242 (1982)
\endref

\ref\no{[P-R]}\by V.~Platonov and A.~Rapinchuk
\jour ``Algebraic groups and Number Theory'', A.P. 1991
\endref

\ref\no\by M.~Queffeler\no{[Quef]}
\paper
Substitution Dynamical systems 
\jour Spectral Analysis, LNM 1294
\endref

\ref\no{[R1]} \by M.~Ratner
\jour Annals Vol 118, 277--313, 1983
\endref

\ref\no{[R2]} \by M.~Ratner
\jour Annals 134, 545--607 (1991)
\endref

\ref\no{[Sa1]}\by P.~Sarnak
\jour Three lectures on the Mobius Function randomness and dynamics, publications.ias.edu/sarnak/
\endref

\ref\no{[Sa2]} \by P.~Sarnak
\jour Clay Math. Proc. Vol 4 (2005), 659--685
\endref


\ref\no {[S-U]}\by P.~Sarnak and A.~Ubis
\paper ``Horocycle flows at prime times''
\jour preprint 2011
\endref

\ref\no{[Va]} \by R.~Vaughan
\jour C.R. Acad. Sci A 285, 981--983 (1977)
\endref

\ref\no{[Ve1]} \by W.Veech
\paper Gauss measures for transformations on the space of interval exchange maps
\jour Annals Vol 115, 210--242 (1982)
\endref

\ref\no{[Ve2]} \by W.Veech
\paper The metric theory of interval exchange transforms
\jour I, II, III, Amer. J.~Math. 106 (1984), 1331-1421
\endref

\ref\no {[V]} \by I.M.~Vinogradov
\jour Recueiv Math (1937)
\endref

\ref\no{[We]} \by A.~Weil
\jour J.~Ind.~Math. Soc. 24 (1961), 589--623
\endref
\endRefs
\enddocument